\documentclass[reqno,12pt]{amsart}
\usepackage{graphics}
\usepackage{a4,amsmath,epsfig}

\def\goth#1{{\mathfrak #1}}

\def\0{\goth 0}
\def\gA{\goth A}
\def\gB{\goth B}

\begin{document}

\vspace*{0.2cm}

\begin{center}
{\Large \bf On models of non-Eucledian spaces generated by associative algebras}\\[1.5 cm]

{\bf\large Maria Trnkov\'a$^{1,2}$}\\[0.4cm] email: {\tt
M.D.Trnkova@gmail.com}\\[0.3cm]
{\it$^1$Department of Algebra and Geometry, Faculty of Science,}\\
{\it Palacky University in Olomouc, Tomkova 10, Olomouc, Czech
Republic}\\[0.1cm]
{\it$^2$Department of Geometry, NIIMM,}\\
{\it Kazan State University, prof. Nuzhina 1/37, Kazan, Russia}
\\[0.5 cm]

\end{center}

\vspace*{1.0cm}

\textbf{Abstract} \emph{ We present the non-trivial example how to
generate non-Euclidean geometries from associative unital algebras.
We consider bundles of the sphere of the degenerate non-Eucleadian
space and its two models. The first (conformal) model is obtained by
the mapping $S$ onto a plane pass through the origin. It is
analogous to the stereographic mapping. The second model
(projective) is constructed by the Norden normalization method,
where we project the sphere onto a plane of normalization defining
the metric and Christoffel symbols which allow us to find geodesic
curves.}

\section{Introduction}

A lot of models of non-Eucledian spaces were studied in the past,
especially spaces of a constant curvature, projective spaces and the
conformal planes (e.g. \cite{Redei}, \cite{Ros}, \cite{RosI},
\cite{Yag}). There exists a lot of studies on how these models can
be generated by algebras. It is well known that algebras define some
structures in bundle manifolds of different types (e.g. \cite{Shir},
\cite{Pavlov}, \cite{Malakh}). In the literature, we can find many
applications of this approach on the cases of non-Eucledian spaces
(e.g. \cite{Vish}, \cite{Kiri}, \cite{Morim}, \cite{Shur},
\cite{Yano}).

We would like to present non-standard models within this framework.
In the beginning, we describe how an associative algebra generates a
vector space and we also discuss some of its properties. In the next
section we define a sphere and the map $S$ in this vector space and
we use it to construct a conformal model. In the third section we
remind some facts about the Norden normalization method \cite{Nord}
and we use it for the construction of a projective model.

Let us denote by $\gA$ a unital associative $n$-dimensional algebra
with the multiplication $xy$, and by $G\subset\gA$ the set of
invertible elements. Then  $G$ is a Lie group with the same
multiplication rule. Let $\gB\subset \gA$ be a unital subalgebra of
an algebra $\gA$ and $H\subset G$ be the set of its invertible
elements. So, $H$ is a Lie subgroup of group $G$ and $G/H$ is the
factor-space of right cosets. A bundle $(G, \pi, M=G/H)$ is a
principal bundle with the structure group $H$, where $\pi$ is a
canonical projection (for example, \cite{Bour}, \cite{Hus}).

Foundations of the theory of finite-dimensional associative algebras
were made by E.~Cartan (1898), Wedderburn (1908) and F.~E.~Molin
(1983), who discovered the structure of any algebra over an
arbitrary base field \cite{Encyc}. E.~Study and E.~Cartan in
\cite{Stud} classified all 3 and 4-dimensional unital associative
irreducible\footnote{{\it Irreducible} means indecomposable into a
direct sum of algebras.} algebras up to an isomorphism. This
classification could be also found in \cite{VishI}. In this paper we
consider only one type of 3-dimensional algebra $\gA$. We leave a
more complicated 4-dimensional case for a future work.

Let $\{\mathbf{1},e_1,e_2\}$ be a basis of our algebra $\gA$ with
the identity element $\mathbf{1}$. The multiplication rules are:
\begin{equation}\label{multiplication}
(e_1)^2=\mathbf{1},\; (e_2)^2=0,\; e_1e_2=-e_2e_1=e_2.
\end{equation}
The matrix representation of an algebra $\gA$ is a space of upper
triangular matrices $T_u=\{\left(\begin{array}{cc}
             x_0 & x_2 \\
             0 & x_1 \\
             \end{array}
       \right)
|x=x_0+x_1\cdot e_1+x_2\cdot e_2\in \gA\}$ with the basic elements
\cite{Encyc}
\begin{equation}\label{basis}
\mathbf{1}=\left(
                       \begin{array}{cc}
                         1 & 0 \\
                         0 & 1 \\
                       \end{array}
                     \right),\quad
                     e_1=\left(
                       \begin{array}{cc}
                         1 & 0 \\
                         0 & -1 \\
                       \end{array}
                     \right),\quad
                     e_2=\left(
                       \begin{array}{cc}
                         0 & 1 \\
                         0 & 0 \\
                       \end{array}
                     \right).
\end{equation}

We consider the trivial conjugation $x=x_{0}+x^{i}e_{i}\to
\overline{x}=x_{0}-x^{i}e_{i}$ with the property
$\overline{xy}=\overline{y}\,\overline{x}$ and the bilinear form

\begin{equation}\label{bilinear form}
(x,y)=\frac{1}{2}(x\overline{y}+y\overline{x}).
\end{equation}

This form takes the real values and it determines a degenerate
scalar product:

\begin{equation}\label{scalar product}
(x,y)=x_{0}y_{0}-x_{1}y_{1}.
\end{equation}

It defines a structure of degenerated pseudo-Euclidean vector spaces
with rank 2 in the algebra $\gA$. (It is also possible to call this
space as "semi-pseudo-Euclidean", but later we will call it just
"pseudo-Euclidean".) The set of invertible elements $G=\{x\in
{\gA}\mid (x_{0})^{2}-(x_{1})^{2}\ne 0\}$ is a non-Abelian Lie
group. Its underlying manifold is $\mathbb{R}^{3}$ without two
transversal 2-planes, hence it consists of 4 connected components.

The norm is defined as usual, $|x,y|^2 = (x-y,x-y)$. The geodesic
curves $x(t)$ are then
$$
x_0 = a_0 t + b_0\qquad x_1 = a_1 t + b_1 \qquad x_2 = f(t)
$$
where $f(t)$ is an arbitrary function of $t$ and $a_0$, $a_1$,
$b_0$, $b_1$ are the numerical coefficients.

In the basis (\ref{basis}) we can find two subalgebras: $R(e_1)$
with basis $\{\mathbf{1}, e_1\}$, it is an algebra of double
numbers, and a subalgebra $R(e_2)$ with basis $\{\mathbf{1}, e_2\}$,
it is an algebra of dual numbers. The set of their invertible
elements $H_{1}=\{x_{0}+x_{1}e_{1}\in R(e_{1})\mid
x^{2}_{0}-x^{2}_{1}\ne 0\}$ and $H_{2}=\{x_{0}+x_{2}e_{2}\in
R(e_{2})\mid x_{0}\ne 0\}$ are Lie subgroups of the Lie group $G$.

The space of right cosets $H_1x$ defines a trivial principal bundle
$(G, \pi, M=G/H_1)$ over the real line $\mathbb{R}$. The fiber is a
plane without two transversal lines and the structure group is
$H_1$. The manifold of the group $G$ is diffeomorphic to direct sum
$\mathbb{R}\times H_{1}$. The coordinate view of the canonical
projection $\pi$ is:
\begin{equation}\label{pi1} \pi(x)=\frac{x_{2}}{x_{0}-x_{1}}.
\end{equation}
The equation of fibers is:
\begin{equation}\label{fiber1} u(x_{0}-x_{1})-x_{2}=0, \; u\in \mathbb{R}.
\end{equation}

Let us investigate $\mathbb{G}$, the group of transformations of Lie
group $G$. We can easily find that it has no dilations and
inversions while there is a vertical translation $x\rightarrow x+a$,
$a\in G$. Furthermore, $\mathbb{G}$ includes the rotations, resp.
anti-rotations,
$$
x'=ax \qquad \mbox{or} \qquad x'=xa
$$
with $|a|^2=1$, resp. $|a|^2=-1$. These elements can be represented
as:
$$
a=\cosh\varphi \pm \sinh\varphi\, e_1 + u \sinh\varphi\, e_2, \quad
\mbox{resp. } \quad a=\sinh\varphi \pm \cosh\varphi e_1 + u
\cosh\varphi\, e_2,
$$
where $u\in \mathbb{R}$. The anti-rotations transform the elements
with the positive norms into the elements with the negative norms
and visa versa.

The bilinear form (\ref{bilinear form}) in the algebra $\gA$ takes
the real values, therefore it is possible to present it as:
$(x,y)=\frac{1}{2}(x\overline{y}+y\overline{x})=\frac{1}{2}(\overline{x}y+\overline{y}x).$
Consequently, in the case of rotations the hyperbolic cosine of an
angle between $x$ and $x'$ are equal to
\begin{equation}\label{cosh}
\cosh(x,x')=\frac{(x,ax)}{|x||ax|}=\frac{1/2(x\overline{ax}+ax\overline{x})}{|x|^2}=
\frac{1/2(x\overline{x}\,\overline{a}+ax\overline{x})}{|x|^2}=
\frac{1}{2}(\overline{a}+a)=\cosh\varphi,
\end{equation}
and the same for the right multiplication. Similarly we get
$\sinh\varphi$ for anti-rotations. Note that the angle $\varphi$
does not depend on $x$.

Transformations
\begin{equation}\label{rot1}
x'=axb,
\end{equation}
where $|a|^2=\pm1,\; |b|^2=\pm1$, are compositions of rotations
and/or anti-rotations $x'=ax$ and $x'=xb$. We see that (\ref{rot1})
defines {\it proper} rotations and anti-rotations.

Similarly,
\begin{equation}\label{rot2}
x'=a\overline{x}b
\end{equation}
are compositions of the reflection $x'=\overline{x}$ and
transformations (\ref{rot1}). These are {\it improper} rotations and
anti-rotations.

\textbf{Lemma} \emph{Any proper or improper rotation/anti-rotation
of the pseudo-Euclidean space $G$ can be represented by (\ref{rot1})
or (\ref{rot2}).}

\textbf{Proof} Rotations and anti-rotations (\ref{rot1}),
(\ref{rot2}) are compositions of odd and even numbers of reflections
of planes passing through the origin. To each plane corresponds its
orthonormal vector $n$. If vectors $x_1$ and $n$ are collinear, then
$\overline{x}_1n=\overline{n}x_1$ and
$x'_1=-n\overline{x}_1n=-n\overline{n}x_1=-x_1$. If vectors $x_2$
and $n$ are orthogonal, then $\overline{x}_2n+\overline{n}x_2=0$ and
$x'_2=-n\overline{x}_2n=n\overline{n}x_2=x_2$. On the other hand,
any vector $x$ can be represented as a sum of vectors $x_1$ and
$x_2$. It means, that a reflection of the plane is:
$x'=-n\overline{x}n$. Therefore, the composition of even, resp. odd
number of reflections of planes are transformation (\ref{rot1}),
resp. (\ref{rot2}). $\Box$

Translations and rotations/anti-rotations are then isometries. All
transformations can be written in a known form (for further
discussion see e.g. \cite{Yag})
\begin{equation}\label{transf}
\left\{
  \begin{array}{ll}
    \hbox{$x_0'=x_0\cosh\varphi + x_1\sinh\varphi + a_0$} \\
    \hbox{$x_1'=x_1\cosh\varphi + x_0\sinh\varphi + a_1$} \\
    \hbox{$x_2'= u_0x_0 + u_1x_1 + u_2x_2 + a_2$}
  \end{array}
\right.
\\
\end{equation}
where $a=a_ie_i\in G$ and $u_i\in{\mathbb{R}}$.

Let us introduce adapted coordinates $(u,\lambda,\varphi)$ of the
bundle in semi-Euclidean space, here $u$ is a basic coordinate,
$\lambda,\varphi$ are fiber coordinates. If $|x|^{2}>0$, we denote
$\lambda =\pm\sqrt{x_{0}^{2}-x_{1}^{2}}\ne 0$, the sign of $\lambda$
is equal to the sign of $x_0$. The adapted coordinates of the bundle
in this case are:
\begin{equation}\label{adapt>0}
x_{0}=\lambda \cosh \varphi ,\quad x_{1}=\lambda \sinh \varphi ,
\quad x_{2}=u\lambda \exp \varphi ,
\end{equation}
where $\lambda \in\mathbb{R}_{0},\quad u,\varphi\in\mathbb{R}.$

If $|x|^{2}<0$, then we write $\lambda
=\pm\sqrt{x_{1}^{2}-x_{0}^{2}}$, the sign of $\lambda$ is equal to
the sign of $x_1$:
\begin{equation}\label{adapt<0}
x_{0}=\lambda \sinh \varphi ,\quad x_{1}=\lambda\cosh\varphi ,\quad
x_{2}=u\lambda \exp \varphi .
\end{equation}

The structure group acts as follows:
\begin{equation}\label{str-gr1}u'=u, \quad
  \lambda'=\lambda \rho, \quad
  \varphi'=\varphi+\psi,
\end{equation}
where the element $a(0, \rho, \psi)$ of the structure group acts on
the element $x(u, \lambda, \varphi)\in G$. This group consists of 4
connected components.

\section{Conformal model of a sphere}

We call \emph{semi-Euclidean sphere with an unit radius} the set of
all elements of algebra $\gA$ whose square is equal to one,
$$
S^{2}(1)=\{x\in\gA\mid x_{0}^{2}-x_{1}^{2}=1\}.
$$
Analogously, the set of elements with an imaginary unit module
$|x|^{2}=-1$ we call \emph{semi-Euclidean sphere with an imaginary
unit radius} $S^{2}(-1)$. One of these spheres can be obtained from
another one by the rotation.

The transformations (\ref{transf}) are now constrained by additional
relation $x_0^2-y_0^2=1$, therefore, only rotations and vertical
translations remain, $a_0=a_1=0$.

We consider the subbundle of the bundle $(G,\pi,M=G/H_1)$ of
semi-Euclidean sphere $S^{2}(1)$, i.e. the bundle $\pi: S^{2}(1)\to
M$. The fibers of the new bundle are intersections of $S^2(1)$ and
planes (\ref{fiber1}). The restriction of the group of double
numbers $H_1$ to $S^{2}(1)$ is a Lie subgroup $S_1$ of double
numbers with an unit module
$$ S_{1}=\{a_{0}+a_{1}e_{1}\in H_{1}\mid a_{0}^{2}-a_{1}^{2}=1\}\,.
$$ This group consists of two connected components.
The bundle $(S^{2}(1),\pi,M)$ is a principal bundle of the group
$S^{2}(1)$ by the Lie subgroup $S_1$ to right cosets.

We define coordinates adapted to the bundle on semi-Euclidean sphere
$S^{2}(1)$. If $x\in S^{2}(1)$ then from (\ref{adapt>0}) we get
$\lambda =\varepsilon,\, \varepsilon =\pm 1$. The parametric
equation of semi-Euclidean sphere in the adapted coordinates
$(u,\varphi)$ is:
\begin{equation}\label{adapt1}
{\textbf{r}}(u,\varphi)=\varepsilon(\cosh \varphi, \sinh
\varphi,u\exp\varphi),
\end{equation}
where $u$ is a basis coordinate, $\varphi$ is a fiber coordinate.
Different values of $\varepsilon$ correspond to different connected
components of semi-Euclidean sphere $S^{2}(1)$.

Let us define the action of the structure group $S_{1}$ on
semi-Euclidean sphere. From (\ref{str-gr1}) and using the adapted
coordinates of elements $a(0, \varepsilon_{1}, \psi),\, x(u,
\varepsilon, \varphi)\in S^{2}(1)$ we get:
$$u'=u, \quad
  \varepsilon'=\varepsilon \varepsilon_{1}, \qquad
  \varphi'=\varphi+\psi. \qquad
$$
This group also consists of two connected components.

The metric tensor for semi-Euclidean sphere has the matrix
representation:
$$
(g_{ij})=\left(\begin{array}{cc}
0 & 0 \\
0 & -1
\end{array}
\right).
$$
The linear element of the metric is:
\begin{equation}\label{metric1}
ds^{2}_{1}=-d\varphi^{2}.
\end{equation}

Now, we want to define the conformal model of the bundle
$(S^{2}(1),\pi,\mathbb{R})$. For that we need to introduce the
conformal map of the sphere to a disconnected plane
$f:S^{2}(1)\rightarrow Q\in\mathbb{R}^{2}$. $Q$ is located at
$x_0=0$. We know that the sphere consists of two disconnected
components, one with $x_0>0$, and other with $x_0<0$. We choose a
pole at the first one, $N(1,0,0)$. All points of $S^{2}(1)$ except
the line through the pole $N$ are stereographically projected to $Q$
such that the first component of the sphere with $x_0>0$ is mapped
on $x_1=(-\infty,-1)\cup (1,\infty)$ while the second component with
$x_0<0$ is mapped on the strip $x_1=(-1,1)$. We denote $x$, $y$
coordinates on $Q$ such that the $x$ axis lies along $x_1$ while the
$y$ axis along $x_2$. Then
\begin{equation}
\label{stereo1}x=\frac{x_{1}}{1-x_{0}},\quad
y=\frac{x_{2}}{1-x_{0}}\,,
\end{equation}
An inverse map $f^{-1}:\mathbb{R}^{2}\rightarrow S^{2}(1)$ where
$x\neq \pm 1$ is:
\begin{equation}
\label{rev-stereo1}x_{0}=-\frac{1+x^{2}}{1-x^{2}},\quad
x_{1}=\frac{2x}{1-x^{2}},\quad x_{2}=\frac{2y}{1-x^{2}}\,.
\end{equation}
If we substitute formulas (\ref{stereo1}) into (\ref{adapt1}) then
we obtain the relations between coordinates $x, y$ and adapted
coordinates $u,\varphi$ which are on semi-Euclidean sphere:
$$f:\quad x=\frac{\sinh
\varphi}{\varepsilon-\cosh \varphi},\quad
y=\frac{u\exp\varphi}{\varepsilon-\cosh \varphi}\,.$$ Then the
inverse map is:
\begin{equation}
\label{ad-stereo1}
\varphi=\ln\Bigl(\varepsilon\frac{x-1}{x+1}\Bigr),\quad
u=-\frac{2y}{(1-x)^{2}}\,.
\end{equation}

Note that the lines $x=\pm 1$ are not included in the mapping and
$Q$ consists of three disconnected components. Also, the line
$x_0=1, x_1=0$ has no image in this mapping. We add it by hand,
identifying the image of this line with the points $x=\pm\infty$ on
$Q$. Then two disconnected parts $x=(-\infty,-1)$ and $x=(1,\infty)$
are connected and we call this plane $C^2$.

In particular, after enlarging $Q$ into $C^2$ by the infinitely
distant point and ideal line crossing this point, then the
stereographic map $f$ becomes diffeomorphism $S$. Note that the
infinitely distant point is the image of point $N$. The ideal line
is the image of the straight line belonging to $S^{2}(1)$ and
crossing the pole: $x_{0}=1,\, x_{1}=0$.

Let us now consider the commutative diagram:
$$
\begin{array}{rcccccl}
S^{2}(1) & & \stackrel{S}\longrightarrow & & C^{2}\\
& \pi \searrow  & &  \swarrow p \\
 & & \mathbb{R} & & &
\end{array}
$$
The map $p=\pi\circ S^{-1}:C^{2}\rightarrow R$ is defined by this
diagram. We find the coordinate form of this map:
$$u=-\frac{2y}{(1-x)^{2}}\,.$$
The map $p:C^{2}\rightarrow \mathbb{R}$ defines the trivial
principal bundle with the base $\mathbb{R}$ and the structure group
$S_1$.

\textbf{Theorem} \emph{Let $S$ is the map $: S^{2}(1) \rightarrow
C^{2}$ as described before. Then $S$ is a conformal map.}

\textbf{Proof} The metric on $G$ induces the metric on $C^{2}$. In
the coordinates $x, y$ it has the form:
\begin{equation}
\label{metric1'}d\widetilde{s}^{2}=-dx^{2}.
\end{equation}

Let us find the metric of semi-Euclidean sphere from the metric on
$C^{2}$. From (\ref{ad-stereo1}) we get $d\varphi=\frac{2}{x^2-1}dx$
and using (\ref{metric1}) and (\ref{metric1'}) we find:
$$
ds^{2}_{1}=\frac{4}{(x^2-1)^2}d\widetilde{s}^{2}.
$$
Hence, the linear element of semi-Euclidean sphere differs from the
linear element of $C^2$ by a conformal factor and therefore, the map
$S$ is conformal. $\Box$

We find the equation of fibers on $C^{2}$. The 1-parametric fibers
family of the bundle $(S^{2}(1),\pi,\mathbb{R})$ in the adaptive
coordinates (\ref{adapt1}) is: $u=c,\quad c\in \mathbb{R}$. From
(\ref{ad-stereo1}) we get the image of this family under the map
$S$:
\begin{equation}
y=-c/2\cdot(x-1)^2.
\end{equation}
The $C^2$ plane is also fibred by this 1-parametric family of
parabolas.

\section{The projective conformal model}

Now we construct the projective semi-conformal model of the sphere
$S^{2}(1)$ and the principal bundle on it. We use a normalization
method of A.P.Norden \cite{Nord}, \cite{NordI}. A. P. Shirokov in
his work \cite{ShirI} constructed conformal models of Non-Euclidean
spaces with this method.

In a projective space $P_n$ a hypersurface $X_{n-1}$ as an
absolute is called \emph{normalized} if with every point $Q\in X_{n-1}$ there is associated:\\
1) a line $P_I$ which has the point $Q$ as the only intersection
with the tangent space $T_{n-1}$,\\
2) a linear space $P_{n-2}$ that belongs to $T_{n-1}$,
but it does not contain the point $Q$.\\
We call them \emph{normals of first and second types}, $P_I$ and
$P_{II}$.

In order to have a polar normalization, $P_I$ and $P_{II}$ must be
polar with respect to the absolute $X_{n-1}$.

We enlarge the semi-Euclidean space $_{2}E^{3}_{1}$ to a projective
space $P^{3}$. Here $_{k}E^{n}_{l}$ denotes a $n$-dimensional
semi-Euclidean space with the metric tensor of rank $k$, and $l$ is
the number of negative inertia index in a quadric form. We consider
homogeneous coordinates $(y_{0}:y_{1}:y_{2}:y_{3})$ in $P^3$, where
$x_{i}=\frac{y_{i}}{y_{3}}\,, i=0,1,2$. Thus $S^{2}(1):\,
x_{0}^{2}-x_{1}^{2}=1$ describes the hyperquadric in $P^{3}$:

\begin{equation}\label{Q}
y_{0}^{2}-y_{1}^{2}-y_{3}^{2}=0\,.
\end{equation}

Here the projective basis $(E_{0},E_{1},E_{2},E_{3})$ is chosen in
the following way. The vertex $E_{0}$ of basis is inside the
hyperquadric. The other vertices $E_{1},E_{2},E_{3}$ are on its
polar plane, $y_{0}=0$. The line $E_0 E_3$ crosses the hyperquadric
at poles $N(1:0:0:1)$, $N'(1:0:0:-1)$. Vertices $E_{1},E_{2}$ lie on
the polar of the line $E_0 E_3$. The vertex of the hyperquadric
coincides with the vertex $E_2$.

The stereographic map of the projective plane $P^2: y_0=0$ to the
hyperquadric (\ref{Q}) from the pole $N(1:0:0:1)$ is shown on the
picture. Let $U(0:y_{1}:y_{2}:y_{3})\in P^{2}$. If $y_3=0$, then the
line $UN$ belongs to the tangent plane $T_N: y_0-y_3=0$ of the
hyperquadric (\ref{Q}) at the point $N$ and in this case the
intersection point of the line $UN$ with the hyperquadric is not
uniquely determined. If $y_{3}\neq 0$, then the intersection point
of the line $UN$ with the hyperquadric is unique. So, we choose the
line $E_{1}E_{2}: y_{3}=0$ as the line at infinity. In the area
$y_{3}\neq 0$ we consider the Cartesian coordinates
$x_{1}=\frac{y_{1}}{y_{3}}, x_{2}=\frac{y_{2}}{y_{3}}\,.$ Then the
plane $\alpha: y_{0}=0, y_{3}\neq 0$ becomes a plane with an affine
structure $A^{2}$. It is possible to introduce the structure of
semi-Euclidean plane $_{1}E^{2}$ with the linear element
\begin{equation}
ds_{0}^{2}=dx_{1}^{2}\,.
\end{equation}
The hyperquadric and the plane $\alpha$ do not intersect or
intersect in two imaginary parallel lines
\begin{equation}\label{q'} x_{1}^{2}=-1\,.
\end{equation}
The restriction of the stereographic projection to the plane
$\alpha$ maps the point $U(0:x_{1}:x_{2}:1)$ into the point $X_1$
\begin{equation}
X_{1}(-1-x_{1}^{2}:2x_{1}:2x_{2}:1-x_{1}^{2})\,.
\end{equation}
So, the Cartesian coordinates $x_i$ can be used as the local
coordinates at the hyperquadric except the point of its intersection
with the tangent plane $T_N$.

\begin{figure}
\begin{center}
\epsfysize=8cm\epsfbox{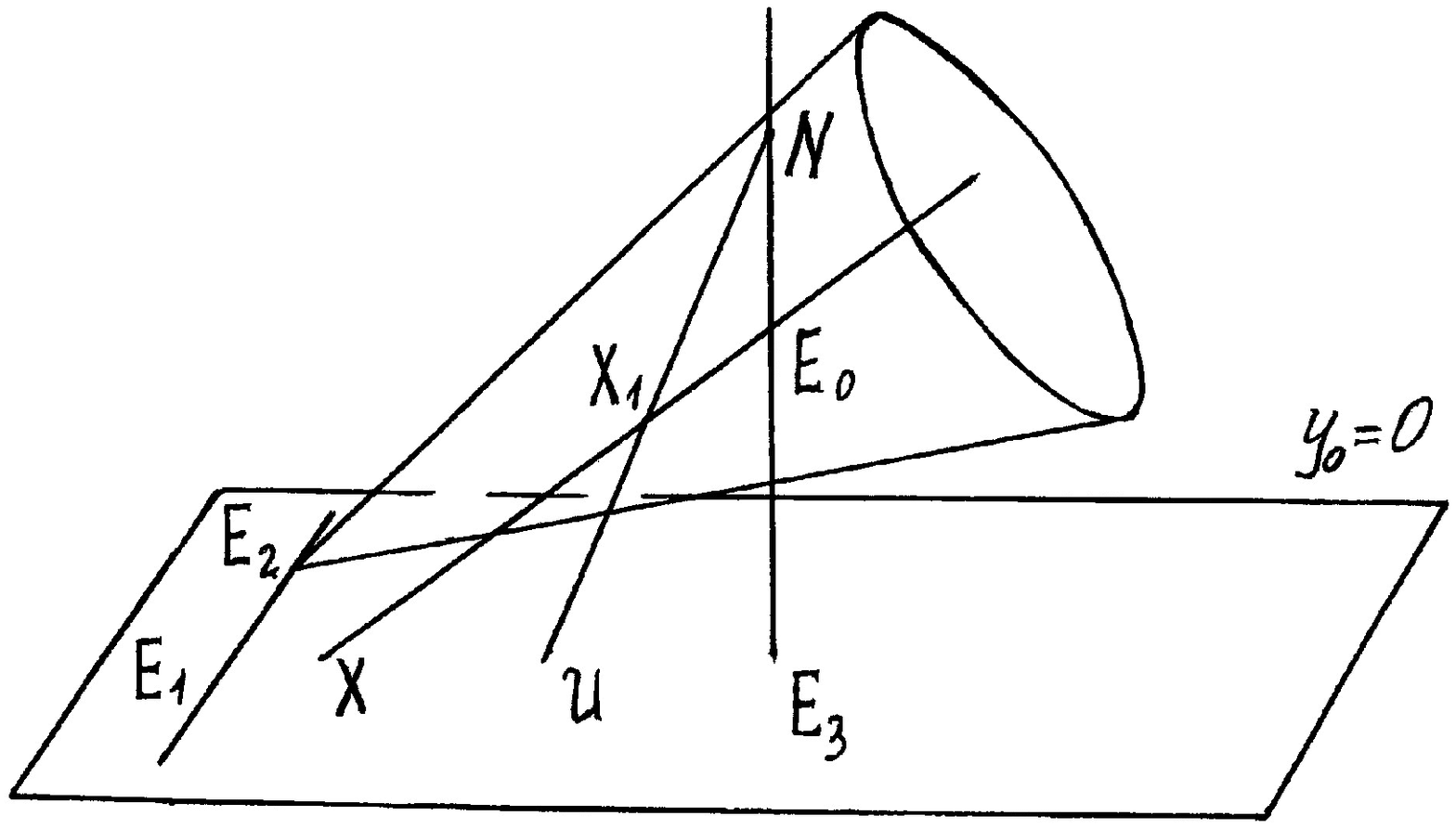}
\end{center}
\end{figure}

We construct an autopolar normalization of the hyperquadric. As a
normal of the first type we take lines with the fixed center $E_0$
and as a normal of the second type we take their polar lines which
belong to the plane $\alpha$ and cross the vertex $E_2$ of the
hyperquadric. The line $E_{0}X_{1}$ intersects the plane $\alpha$ at
the point
$$X(0:2x_{1}:2x_{2}:1-x_{1}^{2})\,.$$ In this normalization the polar
of the point $X$ intersects the plane $\alpha$ on the normal
$P_{II}$. Thus for any point $X$ in the plane $\alpha$ there
corresponds a line which does not cross this point. It means that
the plane $\alpha$ is also normalized. The normalization of $\alpha$
is defined by an absolute quadric (\ref{q'}).

We consider the derivative equations of this normalization. If we
take normals of the first type with fixed center $E_0$, then the
derivative equations (\cite{Nord}, p.204) have the form:
\begin{equation}\label{deriv}
\begin{array}{c}
  \partial_{i}X=Y_{i}+l_{i}X\,, \\
  \nabla_{j}Y_{i}=l_{j}Y_{i}+p_{ji}X\,. \\
\end{array}
\end{equation}
The points $X,\,Y_{i},\,E_{0}$ define a family of projective frames.
Here $Y_{i}$ are generating points of the normal $P_{II}$.

We can calculate the values $(X, X),\; (X,Y_i)$ on the plane
$\alpha$ using the quadric form, which is in the left part of
equation (\ref{Q}). So, $(X,X)=-(1+x_{1}^{2})^2.$

Let us find coordinates of the metric tensor on the plane $\alpha$.
Hence, we take the Weierstrass standardization

$$(\widetilde{X},\widetilde{X})=-1,\quad
\widetilde{X}=\frac{X}{1+x_{1}^{2}}\,.$$

Then the coordinates of the metric tensor are the scalar products of
partial derivatives
$g_{ij}=-(\partial_{i}\widetilde{X},\partial_{j}\widetilde{X})$:
$$
(g_{ij})=\left(%
\begin{array}{cc}
  \frac{4}{(1+x_{1}^{2})^2} & 0 \\
  0 & 0 \\
\end{array}\,.%
\right)
$$
We got the conformal model of the polar normalized plane $\alpha:
y_0=0, y_3\neq0$ with a linear element
\begin{equation}
ds^{2}=\frac{dx_{1}^{2}}{(1+x_{1}^{2})^2}\,.
\end{equation}
It means that this non-Euclidean plane is conformally equivalent to
semi-Euclidean plane $_{1}E^{2}$.

The points $X$ and $Y_{i}$ are conjugated with respect to the polar
(\ref{Q}) and $(X,Y_{i})=0$. From this equation and the derivative
equations (\ref{deriv}) we can get the non-zero connection
coefficients:
$$
\Gamma^{1}_{11}=\Gamma^{2}_{12}=\Gamma^{2}_{21}=-\frac{2x_{1}}{1+x_{1}^{2}}\,,\quad
\Gamma^{2}_{11}=\frac{2x_{2}}{1+x_{1}^{2}}\,.
$$
The sums $\Gamma^s_{ks}=\partial_k \ln \frac{c}{(1+x_1^2)^2}$
($c=const$) are gradients, so the connection is equiaffine.
Curvature tensor has the following non-zero elements:
$$
R_{121\cdot}^{\quad 2}=-R_{211\cdot}^{\quad 2}=-\frac{4}{(1+x_1^2)^2}\,.
$$
Ricci curvature tensor $R_{sk}=R^{\quad i}_{isk\cdot}$ is symmetric:
$R_{11}=\frac{4}{(1+x_1^2)^2}.$ Metric $g_{ij}$ and curvature
$R^{\quad i}_{rsk\cdot}$ tensors are covariantly constant in this
connection: $\nabla_{k}g_{ij}=0,\;\nabla_{l}R^{\quad i}_{rsk\cdot}=
0$. The infinitesimal linear operators for the quadric are
\begin{equation}
\left\{
  \begin{array}{ll}
    \hbox{$L_1 = y_0\frac{\partial}{\partial y_1} + y_1\frac{\partial}{\partial y_0}$} \\
    \hbox{$L_2 = y_0\frac{\partial}{\partial y_3} + y_3\frac{\partial}{\partial y_0}$} \\
    \hbox{$L_3 = y_1\frac{\partial}{\partial y_3} - y_3\frac{\partial}{\partial y_1}$}
  \end{array}
\right.
\\
\end{equation}

Solving geodesic equations we find parametric solutions
\begin{equation}
\left\{
  \begin{array}{ll}
    \hbox{$x_1 = \tan(\omega t + \phi)$} \\
    \hbox{$x_2 = (c_1e^{2i\omega t} + c_2e^{-2i\omega t})\sec^2(\omega t+\phi)$}.
  \end{array}
\right.
\\
\end{equation}
where $c_1$, $c_2$, $\omega$, $\phi$ are integration constants.
Eliminating the parameter $t$ we can rewrite these equations in a
simple form
$$
x_2=A(x_1^2-1)+Bx_1
$$
where $A$ and $B$ are arbitrary constants. We see that the solution
represents parabolas and lines in $x_2x_1$ plane.

Let us consider the bundle of this plane by the double numbers
subalgebra. We write the equations of fibers of semi-Euclidean
sphere $S^2(1)$ in homogeneous coordinates:
\begin{equation}
\left\{%
\begin{array}{ll}
    (y_0-y_1)v-y_{2}=0, \\
    y_{0}^2-y_{1}^2-y_{3}^2=0. \\
\end{array}%
\right.
\end{equation}
This 1-parametric family of curves fibers the hyperquadric and it
defines a bundle on it. The image of these fibers under the
stereographic projection from the pole $N$ to the plane $\alpha$ is:
$$x_2=-v/2\cdot(x_1+1)^2.$$
It is 1-parametric family of parabolas.

\section*{Remark}

We would obtain the similar results for the space of right cosets by
the Lie subgroup $H_2$ (it is the subgroup of invertible dual
numbers) and the bundle of the group $G$ by $H_2$. However, $H_{2}$
is a normal divisor of the group $G$. Therefore, the spaces of right
and left cosets coincide.

\section*{Acknowledgement}
I would like to thank Professor Jiri Vanzura for fruitful
discussions and support with writing this paper. This work was
supported in part by grant No. 201/05/2707 of The Czech Science
Foundation and by the Council of the Czech Government MSM
6198959214.

\end{document}